\documentclass[12pt]{article}

\usepackage{amssymb} 

\newtheorem{theorem}{Theorem}[section]
\newtheorem{lemma}[theorem]{Lemma}
\newtheorem{proposition}[theorem]{Proposition}

\newtheorem{remark}[theorem]{Remark} 
\newtheorem{definition}[theorem]{Definition}

\newcommand{\bull}{\mbox{$\;\;\;$\vrule height .9ex width .8ex depth -.1ex}}
\newcommand{\qed}{$\;\;\;\Box$}
\newenvironment{proof}{\par\smallbreak\noindent{\bf Proof.~}}
{\unskip\nobreak\hfill \bull \par\medbreak}
\newenvironment{proofof}[1]{\par\smallbreak\noindent{\bf Proof of~#1.~}}
{\unskip\nobreak\hfill \bull \par\medbreak}
\newenvironment{subproof}{\par\noindent{\it Proof of Claim.~}}%
{\qed \par\smallbreak}
\newcounter{claim}[theorem]
\renewcommand{\theclaim}{\thetheorem.\arabic{claim}}
\newenvironment{claim}{\refstepcounter{claim}%
\par\medskip\par\noindent{\it Claim~\theclaim.~}~\rm}%
{\par\smallskip\par}

\newcommand{\refeq}[1]{(\ref{#1})}

\newcommand{\setdef}[2]{\left\{ \hspace{0.5mm} #1 :
\hspace{0.5mm} #2 \right\}}
\newcommand{\function}[2]{:#1 \rightarrow #2}

\newcommand{\ehr}[2]{\mathsf{Ehr}(#1,#2)}
\newcommand{\ff}[1]{\mathbb{F}(#1)}
\newcommand{\oo}[1]{\mathbb{O}(#1)}
\newcommand{\dual}[1]{\mathit{Dual}\,#1}
\newcommand{\facing}[1]{\mathit{Facing}\,#1}
\newcommand{\lineg}[1]{\mathcal{L}\,(#1)}
\newcommand{\diam}[1]{\mathit{diam}\,(#1)}
\newcommand{\girth}[2]{\mathit{girth}\,(#1,#2)}
\newcommand{\ssmallskip}{\vspace{1mm}}
\newcommand{\case}[2]{\ssmallskip\par\noindent{\it Case #1:\/ #2}}
\newcommand{\subcase}[2]{\ssmallskip\par\noindent{\it Subcase #1:\/ #2}}
\newcommand{\fingame}[2]{\ssmallskip\par\noindent{\it Endgame #1:\/ #2}}
\newcommand{\loc}[2]{\mathit{loc}_{#1}(#2)}
\newcommand{\glo}[1]{\mathit{glb}(#1)}
\newcommand{\ch}[1]{\mathit{Ch}(#1)}

\title{
On the logical complexity\\ of convex polygon dissections}

\author{Manuel Bodirsky,\ \ \ Mihyun Kang%
\thanks{Supported by DFG Pr 296/7-3.},\ \ \ Oleg Verbitsky%
\thanks{Supported by an Alexander von Humboldt fellowship.}\\[4mm]
Institut f\"ur Informatik,\\
Humboldt Universit\"at zu Berlin, D-10099 Berlin}

\date{21 July 2006}

\begin{document}

\sloppy

\maketitle

\begin{abstract}
The \emph{logical depth} of a graph $G$ is
the minimum quantifier depth of a first order sentence
defining $G$ up to isomorphism in the language of the adjacency and
the equality relations.
We consider the case that $G$ is a dissection of a convex polygon or,
equivalently, a biconnected outerplanar graph. 
We bound the logical depth of a such $G$ from above by a function of
combinatorial parameters of the dual tree of~$G$.
\end{abstract}

\section{Introduction}

A recent series of papers \cite{BFL+,GVe,KPSV,PSV,PSV2,PVV,SSJ,Ver}
investigates a new graph invariant measuring
the complexity of defining a graph in first order logic.
Let $\Phi$ be a first order sentence involving only two relational
symbols, one for vertex adjacency and the other for equality.
The \emph{quantifier depth} (or \emph{rank})
of $\Phi$ is the maximum number of nested quantifiers in it.
We say that $\Phi$ \emph{defines} a graph $G$ if $\Phi$ is true
exactly on those graphs isomorphic to $G$.
The \emph{logical depth} of $G$ is the minimum quantifier
depth of a such $\Phi$. We denote this graph invariant by $D(G)$. 

Estimation of the logical depth of graphs is a natural problem in the 
scope of finite model theory, complexity theory, and graph theory.
There is a close relation of this subject to the graph isomorphism
problem (see, e.g.\ \cite{GVe}).
Many interesting results are obtained even in the particular case
that $G$ is a tree \cite{BFL+,GVe,PSV,Ver}. In \cite{Ver} also
the case of biconnected outerplanar graphs is considered and it is proved
that $D(G)=O(\log n)$ for every $G$ of order $n$ in this class.
The bound is tight up to a constant factor since the cycle of length $n$
has logical depth $\log_2n+O(1)$.

Here we focus on the same class of graphs and prove a general combinatorial
bound for $D(G)$, which implies the worst-case upper bound of \cite{Ver} and,
hopefully, is applicable to estimating $D(G)$ for a random biconnected outerplanar~$G$.

Focusing on the biconnected outerplanar graphs is of interest in several respects.
First, this class admits a natural geometric definition.
Namely, a biconnected outerplanar graph of order $n$ is exactly
a dissection of the $n$-vertex convex polygon viewed as a graph. 
Next, the class of outerplanar graphs
is commonly used as a benchmark for many important algorithmic
problems (see \cite{BKa,BFo,CDR} just as a few examples).
Finally, it occupies an intermediate position between the class of all trees and
the class of
all planar graphs. It is known that the first order 0-1 law is true for the 
former \cite{McC1,McC2} and false for the latter class (see \cite{GNo}
where the probability that a random planar graph has an isolated vertex
is separated both from 0 and 1). 
In particular, this motivates the question whether the 0-1 law is obeyed
for the biconnected outerplanar graphs. Note that such a question
is closely related to the behavior of the logical depth on average~\cite{Spe,KPSV}.

Given a biconnected outerplanar graph $G$,
we estimate $D(G)$ from above by a function of two graph parameters
of the dual tree of $G$, which will be denoted by $\dual G$.
One is the maximum vertex degree, which is as usual denoted by $\Delta(H)$
for a graph $H$. 
Note that $\Delta(\dual G)$ is equal to the maximum length of a facial cycle
in $G$ excluding the outer one.
The other graph parameter, which we borrow from \cite{BFL+}, 
is called the \emph{fineness} of a graph $H$ and denoted by $r(H)$.
Namely, a graph $H$ is called \emph{$r$-fine} if every two vertex disjoint paths of
length $r-1$ in $H$ (on $r$ vertices each) have different degree sequences. 
We define $r(H)$ to be the minimum $r$ such that $H$ is $r$-fine.

\begin{theorem}\label{thm:maiin}
For every biconnected outerplanar graph $G$ on $n$ vertices, we have
$$D(G)\le11\,\log r(\dual G)+5\,\log\Delta(\dual D)+59.$$
\end{theorem}

\noindent
{\bf Organization of the proof and the paper.}
Getting an upper bound for $D(G)$ is equivalent to establishing it
for $D(G,G')$ for all non-isomorphic graphs $G'$, where $D(G,G')$
is the minimum quantifier depth of a sentence distinguishing between
$G$ and $G'$. Dealing with $D(G,G')$ has an advantage that this
number is characterizable as the length of a purely combinatorial
game on $G$ and $G'$, namely, the Ehrenfeucht game, which is a well-known
tool in model theory (see, e.g.\ \cite{Spe}). We define the rules of the game
and state its relations to the first order definability in Section~\ref{s:ehr}.

Assume for a while that $G'$ is also a dissection graph (i.e., biconnected outerplanar).
Estimating $D(G,G')$, it would be natural to try to bound it from above
in terms of $D(\dual G,\dual G')$.
Unfortunately, this is impossible because
non-isomorphic graphs $G$ and $G'$ can have
isomorphic dual trees. By this reason, instead of $\dual G$
we actually deal with a richer structure $\facing G$
that is defined in Section \ref{s:bop} as a structure with two relations
in a special class of structures which we call \emph{graphs with layout}.
Main Lemma 1 proved in Section \ref{s:treelay} provides us with
a general upper bound on $D(T)$ for $T$ a tree with layout in terms of
$\Delta(T)$ and $r(T)$.

Main Lemma 2 stated in Section \ref{s:gvsfacing} bounds $D(G)$ in terms
of $D(\facing G)$, $\Delta(\dual G)$, and $r(\dual G)$.
As one could expect, our overall proof strategy for this lemma consists
in simulating the Ehrenfeucht game on $G$ and $G'$, where
$G$ is a dissection graph and $G'\not\cong G$, by a game on 
$\facing G$ and $\facing{G'}$. An obvious complication is that
$G'$ is not necessarily a dissection graph and therefore we have
either to extend our definition of $\facing{}$ over such $G'$ or to
distinguish $G$ from such $G'$ in another way. We combine both ways. 
In Section \ref{s:gvsfacing} we define a class of pseudo-dissection
graphs and extend the operator $\facing{}$ to this class. We then
distinguish between $G$ and non-pseudo-dissection graph $G'$ in Main Lemma 2A
(Section \ref{s:bopvsnon}), and between $G$ and a non-isomorphic
pseudo-dissection graph $G'$ in Main Lemma 2B (Section \ref{s:bopvspseudo}).
With these preliminaries the proof of Theorem \ref{thm:maiin} in Section \ref{s:uppeer}
takes no efforts. We discuss a potential application of our result in
Section~\ref{s:futur}.

\section{Notation} 

Given a graph $G$, we denote its vertex set by $V(G)$ and edge set by $E(G)$.
The order of $G$ will be denoted by $|G|$, i.e., $|G|=|V(G)|$.
The distance between vertices $u$ and $v$ in $G$ is denoted by $d(u,v)$.
If $u$ and $v$ are in different connected components, we set $d(u,v)=\infty$.
The diameter of $G$ is $\diam G=\max\setdef{d(u,v)}{u,v\in V(G)}$.
The neighborhood of $v$ in $G$, denoted by $\Gamma_G(v)$ 
(the subscript $G$ may be dropped), consists of all vertices
adjacent to $v$. 
The degree of a vertex $v$ is defined by $\deg v=|\Gamma(v)|$.
The maximum degree of $G$ is $\Delta(G)=\max\setdef{\deg v}{v\in V(G)}$.
Given $X\subset V(G)$, we let $G\setminus X$ denote the graph obtained
from $G$ by removal of all vertices in $X$.

$C_n$ stands for the cycle of order $n$. Given a tree and its two vertices
$u$ and $v$, the path from $u$ to $v$ will be denoted by $[u,v]$.
We write $G\cong H$ to say that graphs $G$ and $H$ are isomorphic.

Writing $\log n$, we mean the logarithm base~2.

\section{The Ehrenfeucht game}\label{s:ehr}

Suppose that $G$ and $G'$ are non-isomorphic graphs. 
Let $D(G,G')$ denote the minimum quantifier depth of a first order
sentence $\Phi$ that is true on $G$ but false on $G'$.
The known fact that there are only finitely many pairwise
inequivalent sentences of a fixed quantifier depth easily implies
the equality $D(G)=\max\setdef{D(G,G')}{G'\not\cong G}$.
Thus, estimation of $D(G)$ reduces to estimation of $D(G,G')$.
The latter number has an advantage of having a purely combinatorial
characterization.

In the {\em Ehrenfeucht game\/} on two vertex disjoint graphs $G$ and $G'$
two players, Spoiler and Duplicator, alternatingly select vertices of the graphs,
one vertex per move.
Spoiler starts and is always free to move in any of $G$ and $G'$; Then Duplicator
must choose a vertex in the other graph. Let $x_i\in V(G)$ and $y_i\in V(G')$
denote the vertices selected by the players in the $i$-th round. Duplicator
wins the $k$-round game if the component-wise correspondence between
$x_1,\ldots,x_k$ and $y_1,\ldots,y_k$ is a partial isomorphism from $G$ to $G'$;
Otherwise the winner is Spoiler. The game will be denoted by $\ehr G{G'}$.

By the Ehrenfeucht theorem (e.g.\ \cite{Spe}), $D(G,G')$ 
is equal to the minimum $k$ such that Spoiler has a winning strategy
in the $k$-round Ehrenfeucht game on $G$ and $G'$.

The definitions of $D(G)$, $D(G,G')$, and the Ehrenfeucht game on $G$ and $G'$ 
make a perfect sense for $G$ and $G'$ in any class of structures over
a fixed relational vocabulary $\tau$. In the case of graphs $\tau$ consists
of a single binary relation symbol for the adjacency relation. We will also deal
with $\tau$ containing two binary relations and several unary relations 
$U_1,\ldots,U_m$. The binary relations will be always irreflexive and symmetric
and hence can be viewed as two graphs on the common vertex set. The first graph
will be viewed colored: a vertex $v$ received the color $i$ iff $U_i$ is true
on $v$. An isomorphism between two structures over $\tau$ preserves both
binary relations and all the colors. Let $S$ be a $\tau$-structure on
vertex set $V$ and $W\subseteq V$. Then $S$ induces on $W$ the substructure
$S[W]$ where the binary relations are induced on $W$ in the usual graph-theoretic
sense and each element of $S[W]$ keeps the colors it had in $S$.
A partial isomorphism between two $\tau$-structures is an isomorphism
between their substructures. The Ehrenfeucht theorem holds true in this
general setting.

We also need a generalization of the Ehrenfeucht theorem in a different
direction. We now consider a graph coupled with one of its vertices.
We call two such pairs $(G,v)$ and $(G',v')$ {\em isomorphic\/}
if there is an isomorphism from $G$ to $G'$
taking $v$ to $v'$. We say that a formula $\Phi(x)$ with one free
variable $x$ {\em distinguishes\/} $(G,v)$ from $(G',v')$ if $\Phi(x)$ is true
on $G$ with $x$ assigned the value $v$ and false on $G'$ with $x$ assigned
the value $v'$. Given non-isomorphic $(G,v)$ and $(G',v')$, let
$D(G,v,G',v')$ denote the minimum quantifier depth of a formula $\Phi(x)$
distinguishing $(G,v)$ from $(G',v')$ (it should be stressed that $G$ and $G'$
itself, without the designated vertices, can quite be isomorphic). 
A variant of the Ehrenfeucht theorem says that $D(G,v,G',v')$ is equal
to the minimum $k$ such that Spoiler has a winning strategy
in the $(k+1)$-round Ehrenfeucht game on $G$ and $G'$ in which the vertices
$v$ and $v'$ are selected in the first round (it is convenient to speak
about a $k$-round game with initial configuration $(v,v')$).

We now collect a couple of Spoiler's tricks. Below we suppose that several rounds
of $\ehr G{G'}$ have been played.
Here and throughout the paper we use a harmless slang: Saying that
Spoiler wins in $k$ moves, we mean that he has a strategy allowing him
to win within $k$ next moves whatever Duplicator's strategy is.

\begin{lemma}[Metric Threat 1\label{lem:mt1}]
Let $u,v\in V(G)$, $u',v'\in V(G')$ and suppose that each pair
$u,u'$ and $v,v'$ is selected in a round. Suppose also that $d(u,v)\ne d(u',v')$
and $d(u,v)\ne\infty$
(in particular, it is possible that $d(u',v')=\infty$).
Then Spoiler wins in at most $\lceil\log d(u,v)\rceil$ moves.
\end{lemma}

\begin{proof}
We may assume that $d(u,v) < d(u',v')$. Spoiler uses the \emph{halving strategy}
(see \cite{Spe} for a detailed account). He selects a vertex $w$ on the halfway
between $u$ and $v$. Whatever Duplicator's response $w'$ is, we have either
$d(u,w)<d(u',w')$ or $d(w,v)<d(w',v')$. This allows Spoiler to do the same trick
again. Eventually Spoiler comes to two adjacent vertices while their
counterparts are non-adjacent, which is his win. Each time the original
distance $d=d(u,v)$ reduces to at most $(d+1)/2$, which bounds the number
of rounds.
\end{proof}

The proof of Lemma \ref{lem:mt1} can be easily adapted to get the
following fact (see \cite[Example 2.3.8]{EFl}).

\begin{lemma}\label{lem:cncm}
Let $n\ne m$. Then
$D(C_n,C_m)\le\lceil\log n\rceil+1$.
\end{lemma}

\begin{lemma}\label{lem:defcn}
$D(C_n)<\log n+3$.
\end{lemma}

\begin{proof}
Consider $\ehr{C_n}G$ with $G\not\cong C_n$.
If $G$ has vertices of degree 1 or more than 2, then Spoiler wins
in 4 moves. If $G$ is a vertex-disjoint union of cycles, then
Spoiler first selects two vertices in different components of $G$
and then wins in at most $\log n$ moves by Lemma \ref{lem:mt1}.
If $G$ is a cycle, Spoiler has a fast win according to Lemma \ref{lem:cncm}.
\end{proof}

\section{Basics on biconnected outerplanar graphs}\label{s:bop}

An \emph{outerplanar graph} is a graph embeddable in a plane so that
all vertices lie on the boundary of the outer (i.e.\ unbounded) face.
Speaking on plane drawings of an outerplanar graph, we will always
assume the latter condition. 
A graph is \emph{biconnected} if it cannot be made disconnected by removal of 
a single vertex.
Throughout the paper we abbreviate the term
\emph{Biconnected OuterPlanar graph} as \emph{BOP graph}.
Since the boundary of every face in a biconnected plane graph
is a cycle, the boundary of the outer face of a plane BOP graph
is a Hamiltonian cycle. It is easy to see that there is no other
Hamiltonian cycle. It follows that all drawings of a BOP
graph are equivalent (see \cite[Section 4.3]{Die} for the definitions). 
Referring to properties of a BOP graph we will always assume its unique,
up to an equivalence, drawing.

\begin{definition}[$\ff G$, $\oo G$, BOP graph's accessories]\label{def:bop}\rm
Let $G$ be a BOP graph. The cycle bounding the outer face will be
called \emph{outer}. The cycle bounding any other face
will be called \emph{facial}. An edge of the outer cycle
will be called \emph{outer}; non-outer edges will called \emph{inner}.
$\ff G$ denotes the set of all facial cycles. $\oo G$ denotes the set
of all outer edges.
\end{definition}

Our very first wish on the way to the main result is to somehow
reduce estimation of $D(G)$ to estimation of $D(\dual G)$, where
$\dual G$ is the dual tree of $G$ (a known concept which we define below).
For trees we already have a pretty large library
of Spoiler's strategies developed in \cite{PSV,Ver} and especially
in \cite{BFL+}. A complication with this approach is that there are
non-isomorphic BOP graphs with isomorphic dual trees. We therefore will
represent the embeddability information about $G$ by a more complex
structure.

\begin{definition}[a graph with layout]\label{def:layout}\rm
Let $H$ be a graph with no vertex of degree 2 and no triangle. 
A binary relation $L$ is called a \emph{layout}
for $H$ if the following conditions are met.
\begin{enumerate}
\item
$L$ is defined on $V(H)$.
\item
$L$ is irreflexive and symmetric.
(Thus, $L$ is formally a graph and we will use the graph-theoretic
terminology for it. However, to distinguish between $H$ and $L$,
we will call the former a graph and the latter a relation.)
\item
For every non-leaf vertex $v$ of $H$, $L$ induces a cycle on $\Gamma_H(v)$.
\item\label{i:uv}
Let $u$ and $v$ be adjacent non-leaf vertices of $H$.
Let $v_1,v_2$ be the neighbors of $u$ in $H$ which are the neighbors of
$v$ in $L$ according to the preceding item. Define $u_1,u_2$ symmetrically.
Then $\{v_1,v_2\}$ and $\{u_1,u_2\}$ do not intersect and either
$\{v_1,u_1\}$ and $\{v_2,u_2\}$ are in $L$ or 
$\{v_1,u_2\}$ and $\{v_2,u_1\}$ are in $L$.
\item
No other pairs of vertices are in $L$ except those mentioned in
the preceding two items.
\end{enumerate}
\end{definition}

\begin{definition}[the facing structure, dual tree, crossing relation]%
\label{def:facing}\rm
Let $G$ be a BOP graph. We define a structure $G=\langle H,L\rangle$,
which is a tree $H$ with layout $L$, as follows.

\emph{Vertices of $H$:}
$V(H)=\ff G\cup\oo G$.

\emph{Edges of $H$:}
$H$ has as many edges as $G$; We simultaneously define a one-to-one
correspondence between $E(G)$ and $E(H)$, called the \emph{crossing relation}.
Two $C,D\in\ff G$ are adjacent in $H$ if they share an edge $e\in E(G)$.
We say that $\{C,D\}\in E(H)$ and $e\in E(G)$ \emph{cross each other}.
Furthermore, $C\in\ff G$ and $e\in\oo G$ are adjacent in $H$ if
$e$ belongs to $C$.
We say that $\{C,e\}\in E(H)$ and $e\in E(G)$ \emph{cross each other}.
Finally, there is no edge within $\oo G$ (thus, $\oo G$ consists
of the leaves of $H$). 

It is easy to see that $H$ is a tree. It is known as the \emph{dual tree} of $G$ and
will be sometimes denoted by $\dual G$. The crossing relation agrees
with the standard way of obtaining the dual of a planar graph.
Given an edge $e$ of $G$ or $H$, we will denote its crossing counterpart
by $e^*$.

\emph{Layout $L$:}
Let $u_1,u_2\in\Gamma_H(v)$. We put $\{u_1,u_2\}$ in $L$ iff
$\{v,u_1\}^*$ and $\{v,u_2\}^*$ are adjacent edges of $G$.
Now, let $v,u,v_1,v_2,u_1,u_2\in V(H)$ be as in Item \ref{i:uv}
of Definition \ref{def:layout}. By definition,
we have to put in $L$ either
$\{v_1,u_1\}$ and $\{v_2,u_2\}$ or $\{v_1,u_2\}$ and $\{v_2,u_1\}$.
We put $\{v_1,u_i\}$ iff $\{u,v_1\}^*$ and $\{v,u_i\}^*$ are adjacent in~$G$.
\end{definition}

The structure $\facing G$ for BOP $G$ admits a natural visualization.
Each vertex $C\in\ff G$ should be viewed as a point inside the face of $G$
bounded by $C$. This point sends an edge $e^*$ to (a point inside) every $D\in\ff G$
sharing an edge $e$ with $C$ so that $e^*$ crosses $e$ (note that $C$ and $D$
share at most one edge in a BOP graph). Furthermore, for each outer edge $e$ in $C$,
the (point inside) $C$ sends an edge $e^*$ crossing $e$ and and arriving
at a degree-1 vertex in the outer face. This way we get a plane tree $H$.
The layout $L$ contains all information needed to reconstruct this plane
embedding from $H$ given only as an abstract tree. Namely, for every non-leaf
$v$, $L$ provides the circular order in which the neighbors of $v$
should be arranged in the plane. Given $\Gamma(v)$ already embedded,
$\Gamma(u)$ for an adjacent $u$ could be embedded in two different ways.
$L$ specifies this unambiguously.

What we still have to prove in this section is the following fact
justifying our introducing of the layout notion.

\begin{proposition}\label{prop:facingBOP}  
Let $G$ and $G'$ be BOP graphs. Then $G\cong G'$ if and only if
$\facing G\cong\facing{G'}$, where the latter isomorphism is between
2-relation structures.
\end{proposition}

To prove the proposition it will be helpful to define a kind of 
the line graph of a graph with layout.

\begin{definition}[a chain]\label{def:chain}\rm 
Let $H$ be a graph with layout $L$. Let $v_0,v_1,\ldots,v_k,v_{k+1}$
be a path in $H$. We call it a \emph{chain joining the edges $\{v_0,v_1\}$
and $\{v_k,v_{k+1}\}$} if all $\{v_i,v_{i+2}\}$ and all $\{v_i,v_{i+3}\}$
are in $L$.
\end{definition}

\begin{definition}[the line graph]\rm 
$\lineg G$ denotes the line graph of a graph $G$. 
Given a graph $H$ with layout $L$, we define its \emph{line graph}
$\lineg{H,L}$ in a non-standard way: $\lineg{H,L}$ is a graph whose vertex
set is $E(H)$ and with two $g_1,g_2\in E(H)$ adjacent iff in $H$
there is a chain joining them.
\end{definition}

\begin{proposition}\label{prop:lineiso}  
Let $G$ be a BOP graph. Then $\lineg{\facing G}\cong\lineg G$
and an isomorphism is given by the crossing relation.
\end{proposition}

\begin{proof}
Let $e_1,e_2\in E(G)$. 
We have to show that $e_1$ and $e_2$ have a common vertex $v$ iff
$e_1^*$ and $e_2^*$ are adjacent in $\lineg{\facing G}$.

For the onward direction, let $e_i=\{v,u_i\}$ and $w_1,\ldots,w_{k-1}$
be the neighbors of $v$ lying on the arc of the outer cycle of $G$
between $u_1$ and $u_2$ that does not contain $v$ and listed in the order
from $u_1$ to $u_2$. Let $C_1$ be the facial cycle of $G$ containing 
the edges $e_1$
and $\{v,w_1\}$, $C_i$ the facial cycle containing $\{v,w_{i-1}\}$ and
$\{v,w_{i}\}$, and $C_k$ the facial cycle containing $\{v,w_{k-1}\}$
and $e_2$. Furthermore, let $A$ denote the facial cycle sharing $e_1$ with $C_1$
if such exists or set $A=e_1$ if this is an outer edge.
$B$ is defined similarly for $C_k$ and $e_2$.
It is easy to check that $A,C_1,\ldots,C_k,B$, where $\{A,C_1\}=e_1^*$
and $\{C_k,B\}=e_2^*$, is a chain in $\facing G$.

For the backward direction, suppose that $v_0,v_1,\ldots,v_k,v_{k+1}$ is
a chain in $\facing G$ joining the edges $\{v_0,v_1\}=e_1^*$
and $\{v_k,v_{k+1}\}=e_2^*$. We have to show that $e_1$ and $e_2$ have 
a common vertex $v$. If $k=1$ this follows directly from the 
Definition \ref{def:chain}. If $k=2$ this is also so and, moreover, the
three edges $\{v_0,v_1\}^*$, $\{v_1,v_2\}^*$, $\{v_2,v_3\}^*$ have a common vertex. 
Now let $k=3$. Denote $h_i=\{v_i,v_{i+1}\}^*$. The $h_0,h_1,h_2,h_3$
are pairwise distinct edges of $G$ ($h_0=e_1$ and $h_3=e_2$).
From the case of $k=2$ we already know that $h_0,h_1,h_2$
share a vertex $v$ and $h_1,h_2,h_3$ share a vertex $v'$. As $v$
is the common vertex of $h_1,h_2$ and the same for $v'$, $v=v'$
and we are done. The general case is proved similarly by a simple
induction argument.
\end{proof}

\begin{proofof}{Proposition \ref{prop:facingBOP}}
Let $\phi\function{V(G)}{V(G')}$ be an isomorphism from $G$ to $G'$.
In a natural way $\phi$ extends to a map from $\ff G$ to $\ff{G'}$
and to a map from $\oo G$ to $\oo{G'}$. As it directly follows from
Definition \ref{def:facing}, this is an isomorphism from
$\facing G$ to $\facing{G'}$.

Suppose now that $\facing G\cong\facing{G'}$. Therefore
$\lineg{\facing G}\cong\lineg{\facing{G'}}$ and, by Proposition \ref{prop:lineiso},
$\lineg{G}\cong\lineg{G'}$. Since both $G$ and $G'$ are connected
and neither of them is $K_{1,3}$, the Whitney theorem \cite[Theorem~8.3]{Har}
applies and implies $G\cong G'$.
\end{proofof}

\begin{remark}\label{rem:treetobop}\rm
An easy induction on $n$ shows that every tree $T$ of order $n$ with
layout $L$ uniquely determines a plain BOP $G$ such that
$\langle T,L\rangle=\facing G$. This can be used for an alternative proof 
of Proposition \ref{prop:facingBOP} in the backward direction.
\end{remark}

\section{Complexity of a tree with layout}\label{s:treelay}

\begin{lemma}\label{lem:treelay}
For every colored tree $T$ with layout $L$ we have
$$
D(\langle T,L\rangle) < \log\diam T+\log\Delta(T)+12.
$$
\end{lemma}

\begin{proof}
Let a structure $\langle T',L'\rangle$ consist of
a colored graph $T'$ and a binary relation $L'$ such that 
$\langle T',L'\rangle\not\cong\langle T,L\rangle$.
We design a strategy for Spoiler in $\ehr{\langle T,L\rangle}{\langle T',L'\rangle}$.

\case 1{$T'$ is disconnected.}
Spoiler selects $u'$ and $v'$ in different components of $T'$.
Irrespective of Duplicator's responses $u$ and $v$ in $T$, this poses
Metric Threat 1 and Spoiler wins in no more than $\log d(u,v)+1\le\log\diam T+1$ moves.

\case 2{$T'$ has a cycle.}
Let $C'$ be one of the shortest length. If $|C'|\ge 2\diam T+2$,
Spoiler selects 2 antipodal vertices $u',v'$ of $C'$. As $d(u',v')\ge\diam T+1$,
this is again Metric Threat 1. If $|C'|\le 2\diam T+1$, Spoiler selects
$v'\in C'$ and its neighbors $u',w'$ in $C'$. Duplicator should respond with
a vertex $v$ and its neighbors $u,w$ in $T$. Now Spoiler forces playing
$\ehr{T\setminus v}{T'\setminus v'}$ just never selecting $v$ and $v'$. 
In these graphs $d(u,w)=\infty$ while $d(u',v')\le2\diam T-1$, which again poses
Metric Threat~1.

In the rest of the proof we therefore suppose that $T'$ is a tree.

\case 3{$L'$ is not a layout for $T'$.}
Spoiler is able to reveal any deviation from Definition \ref{def:layout}.
If $T'$ has a vertex of degree 2, Spoiler wins in 4 moves.
If Condition 2 is violated, he needs only 2 moves.
If Condition 3, Spoiler selects a non-leaf $v'$ in $T'$ such that
$L'$ induces on $\Gamma(v')$ a non-cycle. Denote Duplicator's response
in $T$ by $v$. Now Spoiler plays in $\Gamma(v)\cup\Gamma(v')$ and wins
in less than $\log\deg v+3$ moves by Lemma \ref{lem:defcn}.
If anything is wrong with Condition 4, Spoiler wins in 6 moves
(assuming that all the preceding conditions are true).
Finally, if only Condition 5 is false, Spoiler wins in 4 moves.

\case 4{$T'$ is a tree with layout $L'$.}
In our analysis of the game we make two assumptions.

\emph{Assumption 1:} 
Duplicator always replies with a vertex of the same degree. 

This assumption is based on the following claim:
If $v\in V(T)$ and $v'\in V(T')$ are selected in a round and that
$\deg v\ne\deg v'$, then Spoiler wins less than $\log\deg v+2\le\Delta(T)+2$ 
next moves.
Indeed, playing within $\Gamma(v)\cup\Gamma(v')$ and taking into account only
the layouts, he has a fast win due to Lemma~\ref{lem:cncm}.

\emph{Assumption 2:}
Duplicator always respects the graph
metric: For all pairs $u,u'$ and $v,v'$ selected in a round, she takes care
that $d(u,v)=d(u',v')$. (Otherwise she faces Metric Threat 1 and
loses in less than $\log d(u,v)+1\le\log\diam T+1$ moves.)

At the beginning Spoiler selects an arbitrary non-leaf $a\in T$
and two its neighbors $p$ and $q$. Let Duplicator respond with
$a'\in T'$ with $\deg a'=\deg a$ and $p',q'\in\Gamma(a')$. 

While $d$ stands for the graph metric on $T$ and $T'$, we will use
notation $d_L$ for the graph metric on $L$ and $L'$.
Note that the shortest path between $p$ and $q$ in $L$ lie in $\Gamma(a)$
and the similar statement is true for $L'$.
If $d_L(p',q')\ne d_L(p,q)$, this is Metric Threat 1 with respect to 
the layout relation and Spoiler wins in less than $\log d_L(p,q)+1\le\log\deg a$
moves. We hence assume in the sequel that $d_L(p',q')=d_L(p,q)$.

We now introduce some notation
and notions for $T$ which carry through in an obvious way also to $T'$.

The layout $L$ determines a plane embedding of $T$, which is unique up to
a continuous bijective map of the plain onto itself (cf.\ Remark \ref{rem:treetobop}).
Namely, we first place an arbitrary non-leaf of $T$ in the plane and then arrange
its neighborhood according to the circular order given by $L$.
Having embedded a non-leaf $x$ and $\Gamma(x)$, for each non-leaf $y\in\Gamma(x)$
we can embed $\Gamma(y)$ (as a plain $L$-cycle) in two ways. We always choose the one
under which the two $L$-edges between $\Gamma(x)$ and $\Gamma(y)$ do not
intersect one another (or, equivalently, 
they are on the different sides of the edge $\{x,y\}$).

For our convenience we suppose that the direction from $p$ to $q$ around $a$
is counter-clockwise. Fixing this convention justifies using some space-orientation
terminology. For example, let $\deg u\ne1$ and $t\in\Gamma(u)$. The $t$ has two
$L$-neighbors in $\Gamma(u)$. We will call one of these two, denoted by $l_u(t)$, 
\emph{left}
if the direction from $t$ to $l$ around $u$ is counter-clockwise ($l$ lies on the left
from $t$ if looking from $u$). The second will be called the \emph{right
$L$-neighbor of $t$ in $\Gamma(u)$} and denoted by $r_u(t)$.

For each non-leaf $u$ we intend to coordinatize
$\Gamma(u)$ and, for this purpose, we want to designate two neighbors of $u$,
that will be denoted by $p_u$ and $q_u$. 
If $u=a$, we set $p_u=p$ and $q_u=q$. If $u\ne a$, let $t$ precede $u$ on the way
from $a$. Then we set $p_u=t$ and $q_u$ to be the left $L$-neighbor
of $w$ in $\Gamma(u)$. Given $v\in\Gamma(u)$, we define the
\emph{local coordinates of $v$ in $\Gamma(u)$} by 
$\loc  uv=(d_L(v,p_u),d_L(v,p_v))$. 

We define the \emph{global coordinates of $v\in V(T)$}
(with respect to $a,p,q$) recursively dependent on $d(a,v)$.
First, $\glo a$ is empty. If $v\ne a$, let $u$ precede $v$ on the way
from $a$. To obtain the global coordinates of $v$, we concatenate
the global coordinates of $u$ and the local coordinates of $v$ in $\Gamma(u)$:
$\glo v=\glo u\loc uv$.

Since the structures under consideration are non-isomorphic, there exist
$v\in V(T)$ and $v'\in V(T')$ such that $\glo v=\glo{v'}$ but
either $\deg v\ne\deg v'$ or $v$ and $v'$ have different colors.
In the fourth round Spoiler selects such a $v$. Denote Duplicator's
response by $v'$. If $v$ and $v'$ have different colors, she loses immediately.
By Assumption 1, $\deg v=\deg v'$.
We therefore assume that $\glo v\ne\glo{v'}$. 
By Assumption 2, $d(a,v)=d(a',v')$.

Let $w$ be the first vertex on the way from $a$ to $v$ such that
$\glo w\ne\glo{w'}$ for the $w'$ on the path between $a'$ and $v'$ 
with $d(a,w)=d(a',w')$. Denote the predecessors of $w$ and $w'$ by
$u$ and $u'$ respectively. Thus, $\glo u=\glo{u'}$ and 
\begin{equation}\label{eq:locne}
\loc uw\ne\loc{u'}{w'}.
\end{equation}

In the fifth and sixth rounds Spoiler selects $u$ and $w$.
By Assumption 2, 
Duplicator answers with $u'$ and $w'$.
If $u=a$ and $u'=a'$, Inequality \refeq{eq:locne} gives Metric Threat 1
in the layout cycles on $\Gamma(a)$ and $\Gamma(a')$, which allows
Spoiler to win in less than $\log\deg a+1\le\log\Delta(T)+1$ moves.
We hence assume that $u\ne a$, i.e., $d(a,w)>1$.

Denote also the predecessor of $u$ (resp.\ $u'$) on the way from $a$
(resp.\ $a'$) by $t$ (resp.\ $t'$).
In the seventh round Spoiler selects $t$ and,
by Assumption 2, Duplicator should answer with $t'$.
In the eighth round Spoiler selects $z=q_u$. Denote Duplicator's answer by $z'$.

\subcase{4.1}{$z'\notin\{l_{u'}(t'),r_{u'}(t')\}$.}
This is Spoiler's win as the layout relation is violated.

\subcase{4.2}{$z'=l_{u'}(t')$ ($z'=q_{u'}$).}
Note that  $u$, $q_u$, and $p_u=t$, the local coordinate origins in $\Gamma(u)$,
as well as their counterparts $u'$, $q_{u'}=z'$, and $p_{u'}=t'$, the local 
coordinate origins in $\Gamma(u')$, are selected. In view of \refeq{eq:locne} 
this poses Metric Threat 1 and Spoiler wins in less than $\log\Delta(T)+1$ moves.
The only remaining case is as follows.

\subcase{4.3}{$z'=r_{u'}(t')$.}
We first consider the case that $d(a,w)=2$, i.e., $t=a$.
In the ninth round Spoiler selects the vertex $y\in\Gamma(a)$ which is
in the relation $L$ to $z=q_u$ ($y$ is one of the two $L$-neighbors of $u$
in $\Gamma(a)$). Not to lose immediately, Duplicator selects
the vertex $y'\in\Gamma(a')$ which is in the relation $L'$ to $z'$ 
($y'$ is one of the two $L'$-neighbors of $u'$ in $\Gamma(a')$).
Since $\glo u=\glo{u'}$, we have $\loc au=\loc{a'}{u'}$.
Recall that, while $z'$ is the right $L'$-neighbor of $a'$ in $\Gamma(u')$,
$z$ is the left $L$-neighbor of $a$ in $\Gamma(u)$. It follows that
$\loc ay\ne\loc{a'}{y'}$ and therefore Spoiler wins in less than $\log\deg a$
moves.

Suppose that $d(a,w)>2$. Let $b$ (resp.\ $b'$) denote the successor of
$a$ (resp.\ $a'$) on the way to $u$ (resp.\ $u'$). In the ninth round Spoiler 
selects $q_b$. Duplicator should respond with
$q_{b'}$ for else she loses in less than $\log\deg a+2$
moves similarly to the case that $d(a,w)=2$.

Given a vertex $x\ne a$ on the path from $a$ to $u$, we let
$\hat x$ denote the predecessor of $x$ on this path. Furthermore,
we call the $l_x(\hat x)$ a \emph{left vertex} and $r_x(\hat x)$ 
a \emph{right vertex}.
From now on Spoiler follows a kind of the halving strategy.
At the beginning he has two left vertices selected in $T$,
$q_b$ and $q_u$. The corresponding vertices selected in $T'$ are
$q_{b'}$ and $z'$, the former being the left and the latter
being the right vertex. Spoiler is able to make the distance between
such vertices (left-left in $T$ and left-right in $T'$) shorter.
He selects the vertex $s=l_x(\hat x)$ for an $x$ on the halfway between $b$ and $u$ 
(at most two choices of $x$ are possible). Denote Duplicator's response by $s'$.
Let $x'$ denote the analog of $x$ on the path from $b'$ to $u'$.

If $s'$ is not an $L'$-neighbor of $\hat x'$ in $\Gamma(x')$, then,
by Assumption 2, Spoiler wins in the next two moves selecting $\hat x'$ and $x'$.
We hence suppose that $s'=l_{x'}(\hat x')$ or $s'=r_{x'}(\hat x')$.
In either case one of the pairs $q_{b'},s'$ and $s',z'$ consists of
a left and a right vertex, whereas the corresponding pair in $T$
consists of two left vertices. Note that 
$\max\{d(q_b,x),d(x,q_u)\}\le(d(q_b,q_u)+1)/2$.
It follows that, if Spoiler repeats the same trick with the new pairs,
then in less than $\log d(q_b,q_u)+1$ moves he ensures that,
for two adjacent $x_1,x_2$ on the path between $b$ and $u$
and for their analogs $x'_1,x'_2$ on the path between $b'$ and $u'$,
selected are $l_{x_1}(\hat x_1)$ and $l_{x_2}(\hat x_2)$ in $T$
and $l_{x'_1}(\hat x'_1)$ and $r_{x'_2}(\hat x'_2)$ in $T'$.

Assume that $x_1$ is closer to $b$ than $x_2$ is (the other case is similar).
In two more moves Spoiler selects $x_1$ and $\hat x_1$.
Note that the selected vertices $\hat x_1$ and $l_{x_1}(\hat x_1)$
are the local coordinate origins in $\Gamma(x_1)$. This allows Spoiler
to proceed as in the case of $d(a,w)=2$ (with $\Gamma(x_1)$ in place of
$\Gamma(a)$ and $x_2$ in place of $u$) and win in at less than 
$\log\deg x_1+1\le\log\Delta(T)+1$ extra moves.
\end{proof}

\begin{lemma}[Main Lemma 1]\label{lem:ml2}
For every tree $T$ with layout $L$ we have
$$
D(\langle T,L\rangle)\le3\,\log r(T)+\log\Delta(T)+18.
$$
\end{lemma}

In our proof of Lemma \ref{lem:ml2}
we borrow some ideas from the proof of \cite[Theorem 19]{BFL+}. 

\begin{definition}
Given a graph $G$ and an integer $r$, call its vertex $v$ an \emph{$r$-yuppie}
if there are vertices $u$ and $w$ such that $d(u,w)=2r$ and $d(u,v)=d(v,w)=r$.
Let $Y_r(G)$ denote the set of $r$-yuppies in~$G$.
\end{definition}

Note that there is a first order formula $\Phi_r(x)$ of quantifier depth
less than $\log r+4$ such that $v\in Y_r(G)$ iff $\Phi_r(x)$ is true on $G,v$.

Note also that, if $G$ is a tree, then $Y(G)$ spans a subtree.

\begin{lemma}\label{lem:yuppie}
Let $T$ be a tree with layout $L$. 
Let $r=r(T)$.
Suppose that $v,v'\in Y_r(T)$
and $v\ne v'$. Then 
$$
D(\langle T,L\rangle,v,\langle T,L\rangle,v')<\max\{\log r(T),\log\Delta(T)\}+7.
$$
\end{lemma}

\begin{proof}
Given $w,w^*\in V(T)$ with $d(w,w^*)=r$, let $w,w_1,\ldots,w_r=w^*$ be
the path of length $r$ from $w$ to $w^*$. Then we set 
$s(w,w^*)=(\deg w_1,\ldots,w_r)$ and define $\ch w$ to be the set of
$s(w,w^*)$ for all $w^*$ at the distance $r$ from~$w$.

\begin{claim}\label{cl:ch}
If $\ch w\ne\ch{w'}$, then 
$D(\langle T,L\rangle,w,\langle T,L\rangle,w')<\max\{\log r(T),\log\Delta(T)\}+4.$
\end{claim}

\begin{subproof}
Without loss of generality suppose that there is $u$ such that
$s(w,u)\notin\ch{w'}$. Spoiler selects this $u$. Denote Duplicator's
response by $u'$. In the next move Spoiler selects a vertex $t\in[w,u]$
whose degree differs from the degree of the corresponding vertex $t'$ in $[w',u']$.
If Duplicator responds not with $t'$, she faces Metric Threat 1
and loses in less than $\log d(w,u)\le\log r+1$ moves. If she responds with
$t'$, Spoiler plays within $\Gamma(t)\cup\Gamma(t')$ taking into
account only the induced layout cycles of lengths $\deg t$ and $\deg t'$
respectively. By Lemma \ref{lem:cncm}, he needs less than 
$\log\deg t+2\le\log\Delta(T)+2$ to win.
\end{subproof}

On the account of the claim we will assume that $\ch v=\ch{v'}$.
Spoiler selects a vertex $a$ at the distance $r$ from $v$ so that
$v'\notin[v,a]$. Let $a'$ be Duplicator's response. Suppose that
$s(v,a)=s(v',a')$ for else Spoiler wins in less than
$\log\max\{r,\Delta(T)\}+3$ moves as in the proof of Claim \ref{cl:ch}.
This equality is possible only if $a'$ lies strictly between $a$ and $v$
and $v$ strictly between $v'$ and $a'$. In the next move Spoiler selects
the neighbor $b$ of $v$ in $[v,a]$. To avoid Metric Threat 1, Duplicator
must respond with $b'$ being the neighbor of $v'$ in $[v',a']$.

Fix a vertex $z$ maximizing $d(v,z)$ under the condition that
$\ch z=\ch v$ and $v\in[b,z]$ ($z=v'$ is not excluded and hence
$d(v,z)\ge d(v,v')$). Similarly, fix a vertex $z'$ maximizing $d(v',z')$ 
under the condition that $\ch{z'}=\ch{v'}$ and $v'\in[b',z']$
(it is not excluded that $z'=v'$). Note that $d(v,z)\ge d(v,z')>d(v',z')$.
In the third move Spoiler selects $c=z$. Denote Duplicator's response
by $c'$. If $\ch{c'}\ne\ch c$, Spoiler wins in less than
$\max\{\log r,\log\Delta(T)\}+4$ moves by Claim \ref{cl:ch}.
Suppose hence that $\ch{c'}=\ch c$. Note that then $\ch c=\ch v$, which 
is possible only if $d(v,c)\le2r$. If $b'\in[v',c']$, this is
Metric Threat 1 and Spoiler wins less than $\log d(v,c)+1\le\log r+2$ moves.
Otherwise $d(v,c)=d(v,z)>d(v',z')\ge d(v',c')$. This is again
Metric Threat 1 and Spoiler wins in less than $\log r+2$ moves.
\end{proof}

\begin{proofof}{Lemma \ref{lem:ml2}}
We will say that a first order formula $\Psi(x)$ \emph{defines a vertex $v$
in a graph $G$} if $\Psi(x)$ is true on $G,v$ but false on $G,u$
for any other $u\in V(G)$. Let $r=r(T)$ and $Y=Y_r(T)$.
According to Lemma \ref{lem:yuppie}, for every two distinct vertices
$v$ and $v'$ in $Y$ we have a first order formula $\Phi_{vv'}(x)$
true on $T,v$, false on $T,v'$, and of quantifier depth less than
$\max\{\log r(T),\log\Delta(T)\}+7$. Note that the formula
$$
\Phi_v(x)\doteq\bigwedge_{v\ne v'\in Y_r(T)}\Phi_{vv'}(x)\wedge\Phi(x)
$$
defines a vertex $v\in Y$ in $T$ with quantifier depth within the same bound.

We analyze $\ehr{\langle T,L\rangle}{\langle F,M\rangle}$ where $F$
is an arbitrary graph with layout $M$ and the two structures are non-isomorphic.
Given $v\in Y$, we call $v'\in V(F)$ a \emph{friend of $v$} if
$F,v'$ satisfies $\Phi_v(x)$. 

\emph{Assumption.} Duplicator respects the friendship: She takes case
that, if $v\in Y$ and $v'\in V(F)$ are selected in a round,
then they are friends. Once this is not so, Spoiler wins in at most
$D(\Phi_v)$ moves.

Under Assumption, Spoiler wins just in 2 moves in the case that
some $v\in Y$ has no or more than one friends. We therefore
assume that the friendship relation is a one-to-one correspondence
between $Y$ and a set $Y'\subset V(F)$. Actually $Y'\subseteq Y_r(F)$
and, if the inclusion is proper, Spoiler wins fast. Namely, he selects
a vertex $v'\in Y_r(F)\setminus Y'$. Respecting the friendship, 
Duplicator moves outside $Y$ after which Spoiler wins in $D(\Phi)$
moves. We hence suppose that $Y'=Y_r(F)$.

Denote the substructures induced by $\langle T,L\rangle$ and
$\langle F,M\rangle$ on $Y$ and $Y'$
by $\langle K,M\rangle$ and $\langle K',M'\rangle$ respectively.
We will assume that the friendship gives an isomorphism between
$\langle K,M\rangle$ and $\langle K',M'\rangle$ (for else
Spoiler wins in two extra moves).

Given $v\in Y$, let $T_v$ denote the subgraph of $T$
spanned by those $u$ for which there is a path from $v$ to $u$
vertex-disjoint with $Y\setminus\{v\}$. Given $v'\in Y'$, the
subgraph $F_{v'}$ of $\langle F,M\rangle$ is defined similarly.
Note that $T_v$ is a tree.

\case 1{$F$ has a vertex $w'$ that belongs neither to $Y'$ nor to
$F_{v'}$ for any $v'\in Y'$.}
Spoiler selects such a $w'$. Denote Duplicator's answer by $w$.
The $w$, as every vertex of $T$, has the following property: 
there is a yuppie $u$ at the distance at most $r$ from $w'$.
We are done because this property is first order expressible
with quantifier depth less than $\log r+5$.

\case 2{$\diam{F_{v'}}>2r$ for some $v'\in Y'$.}
In this case $F_{v'}$ has a vertex $w'$ at the distance $r+1$ from $v'$.
Spoiler selects $v'$ and $w'$. Denote Duplicator's responses by $v$ and $w$. 
By Assumption, $v$ is the friend of $v'$ in $Y$. We also assume that
$w\notin Y$ for else Spoiler wins in $D(\Phi)$ moves.

If $w\in T_v$, then $d(v,w)\le r$ in $T\setminus Y$ while
$d(v',w')\le r$ in $F\setminus Y'$. Spoiler is here able to implement
Metric Threat 1 and win in less than $\log r+1$ moves because,
by Assumption, Duplicator respects the non-membership in $Y$ and $Y'$.

If $w\notin T_v$, then there is no path in $T\setminus Y$ connecting 
$w$ and $v$ and Spoiler as well is able to implement Metric Threat 1
with $F$ and $T$ interchanged.

\case 3{Not Case 2 and some $F_{v'}$ and $F_{u'}$ intersect.}
Let $w'$ be a common vertex. Spoiler selects $v'$, $u'$, and $w'$.
Let Duplicator respond with $v,u\in Y$ and $w\notin Y$.
As we are not in Case 2, in $F\setminus Y'$ there are paths
from $w'$ to $v'$ and $u'$ of length at most $2r$ each.
This is not so in $T$ for $w$ and $v$ or for $w$ and $u$
and hence Spoiler wins in less than $\log r+2$ moves.

We color the vertices of $Y$ and $Y'$ in the following
way: All vertices in $Y$ have pairwise
distinct colors and the coloring of $Y'$ is determined by the condition
that every two friends have the same color.
For each $v\in Y$, we now define $T^*_v$, a colored graph with layout,
to be the substructure of (the colored version of) $\langle T,L\rangle$
induced on the set $V(T^*_v)=V(T^*_v)\cup\Gamma_T(v)$. 
For each $v'\in Y'$, a colored graph with layout $F^*_{v'}$
is defined similarly.

Recall that $\langle T,L\rangle\not\cong\langle F,M\rangle$.
It is not hard to see that, if we are in none of Cases 1--3,
we must have the following.

\case 4{There is $v\in Y$ such that for its friend $v'\in Y'$
the colored graphs with layouts $T^*_v$ and $F^*_{v'}$ are non-isomorphic
(and none of Cases 1--3).}
Spoiler plays $\ehr{T^*_v}{F^*_{v'}}$. Once Duplicator deviates selecting
a vertex $u\notin T^*_v$ or $u'\notin F^*_{v'}$, she faces Metric 
Threat 1 because there is no path from $v$ to $u$ through $T^*_v$
and no path from $v'$ to $u'$ through $F^*_{v'}$ whereas Spoiler always
has a path from $v$ of length at most $r+1$ and a path from $v'$
of length at most $2r+1$. Suppose Duplicator agrees to play $\ehr{T^*_v}{F^*_{v'}}$.
Then Spoiler applies the strategy of Lemma \ref{lem:treelay}
and wins in less than $\log(r+1)+\log\Delta(T)+12$ moves.
\end{proofof}

\section{Relationship between $\mathbf{D(G)}$ and $\mathbf{D(\facing G)}$}%
\label{s:gvsfacing}

\begin{definition}[$f(G)$, the facial circumference of a graph]\rm
Given a BOP graph $G$, we denote the maximum length of a cycle in $\ff G$
(the set of facial cycles except for the outer one) by $f(G)$.
\end{definition}

\begin{lemma}[Main Lemma 2]\label{lem:ml1}
If $G$ is a BOP graph, then
$$
D(G) < 3\,D(\facing G)+2\,\log f(G)+2\,\log r(\dual G)+5.
$$
\end{lemma}

The proof of Lemma \ref{lem:ml1} will take two subsequent sections.
The rest of this section is devoted to some preliminaries and discussion
of our proof strategy. 

It would be much easier to prove a weaker fact:
if $G$ and $G'$ are two non-isomorphic BOP graphs, then
\begin{equation}\label{eq:ggfacing}
D(G,G')\le 3\,D(\facing G,\facing{G'})+2\,\log f(G)+2\,\log r(\dual G)+O(1).
\end{equation}
However, we have to take into account also that $G'$ may be a non-BOP
graph. It would hence be desirable to have two separate fast strategies
for Spoiler in $\ehr{G,G'}$ for the case that $G'$ is BOP and the case that
$G'$ is not, where the ``fast strategy'' has finally result in a 
$O(\log\log n)$-bound. Unfortunately, this is impossible since, for example,
for the cycle and its two vertex disjoint copies we have
$D(C_n,2C_n)=\log n-O(1)$ \cite[Example 2.3.8]{EFl}.

Our way is therefore longer. First, we will extend the class of BOP graphs to
a broader class of \emph{pseudo-BOP graphs} so that $D(G,G')=O(\log\log n)$
for every BOP graph $G$ and non-pseudo-BOP graph $G'$. Second, we will be
able to extend the operator $\facing{}$ to pseudo-BOP graphs and prove
\refeq{eq:ggfacing} for every BOP graph $G$ and pseudo-BOP graph $G'$.

\begin{definition}[the shortest cycle via two vertices,
a pseudo-facial cycle]\label{def:pseudofl}\rm
Let $u,v\in V(G)$. We define $\girth uv$ to be the minimum length
of a cycle in $G$ going through $u$ and $v$. We call $C$ \emph{the shortest
cycle via $u$ and $v$} if $C$ is a unique cycle of length at most $\girth uv$ 
going through $u$ and~$v$.

We call a cycle $C$ \emph{pseudo-facial} if, for every two vertices $u$ and $v$
in $C$ that are non-adjacent in $C$, $C$ is the shortest cycle via $u$ and~$v$.

The set of pseudo-facial cycles of $G$ will be denoted by $\ff G$.
The maximum length of a cycle in $\ff G$ will be denoted by~$f(G)$.
\end{definition}

\begin{remark}\rm
\mbox{}

\begin{enumerate}
\item
Every triangle is a pseudo-facial cycle.
\item
Every pseudo-facial cycle is chordless.
\item
In a BOP graph the notions of a facial cycle and a pseudo-facial cycle coincide.
Hence the notation $\ff G$ in Definitions \ref{def:bop} and \ref{def:pseudofl}
is coherent.
\end{enumerate}
\end{remark}

\begin{proposition}\label{prop:atmosttwo}
Two pseudo-facial cycles have at most two vertices in common and,
if they have two, those are a common edge.
\end{proposition}

\begin{proof}
We first observe that, if two pseudo-facial cycles have two vertices 
$u$ and $v$ in common, then either of them contains the edge $\{u,v\}$.
Indeed, if one cycle has this edge, the other cycle has it as well
since pseudo-facial cycles are chordless.
The case that $\{u,v\}$ is an edge in neither cycle is impossible
since then $u$ and $v$ would have two cycles of the minimum length,
contradicting Definition~\ref{def:pseudofl}.

It remains to show that two pseudo-facial cycles cannot have 3 common
vertices. This easily follows from the above observation and the fact
that pseudo-facial cycles are chordless.
\end{proof}

It is quite non-obvious if Definition \ref{def:pseudofl} provides us with an 
efficiently (i.e.\ polynomial-time)
verifiable criterion of whether a cycle is pseudo-facial or not.
Nevertheless, such a criterion does exist.

\begin{definition}
{\bf(the shortest biconnection of two vertices, a boundary-like cycle)}\rm
\ Let $C$ be a cycle in a graph $G$ containing vertices $u$ and $v$.
Here and throughout, our default convention will be that $P_1$ and $P_2$
denote the paths on which $u$ and $v$ split $C$ and that $|P_1|\le|P_2|$.
If $|P_1|=|P_2|$, we call $u$ and $v$ \emph{antipodal vertices}.

Suppose that $|P_1|\ge2$. We call $C$ \emph{the shortest biconnection of
$u$ and $v$} if the following is true.
\begin{itemize}
\item
If $u$ and $v$ are antipodal, then $|P_1|=|P_2|=d(u,v)$ and,
except $P_1,P_2$, there is no other path of length $d(u,v)$ between $u$ and $v$.
\item
If $u$ and $v$ are not antipodal, then $P_1$ is a unique path of length
$d(u,v)$ between $u$ and $v$ and in $G\setminus(P_1\setminus\{u,v\})$,
except $P_2$, there is no other path of length $\le|P_2|$ connecting $u$ and $v$
in $G\setminus(P_1\setminus\{u,v\})$.
\end{itemize}

We call a cycle $C$ \emph{boundary-like} if, for every two vertices $u$ and $v$
in $C$ that are non-adjacent in $C$, $C$ is the shortest biconnection of $u$ and~$v$.
\end{definition}

\begin{proposition}\label{lem:pfeqbl}
A cycle $C$ is pseudo-facial iff it is boundary-like.
\end{proposition}

\begin{proof}
Assume that $C$ is pseudo-facial. Let $u$ and $v$ be non-adjacent vertices
on $C$ splitting $C$ in paths $P_1$ and $P_2$ with $|P_1|\le|P_2|$.
We have to show that $C$ is the shortest biconnection of $u$ and $v$.
Let $P$ be a shortest path between $u$ and $v$. Assume for a while that
$P\notin\{P_1,P_2\}$. 

If $P$ has no common vertex with $P_1$ and $P_2$
except $u$ and $v$, then this clearly contradicts the fact that $C$
is a unique shortest cycle via $u$ and $v$.
It follows that $P$ contains vertices $u',v'\notin\{u,v\}$ such that
$u'$ is common with $P_1$, $v'$ is common with $P_2$, but no vertex between
$u'$ and $v'$ is shared with $P_1$ or $P_2$. Since $C$ is chordless,
$u'$ and $v'$ are non-adjacent. Notice now that $C$ is not a unique shortest
cycle via $u'$ and $v'$, a contradiction. Thus, $P\in\{P_1,P_2\}$.

This proves that in $G$ there are at most 2 shortest paths
between $u$ and $v$ and that, if there are two, then those are $P_1$
and $P_2$. In the latter case $C$ satisfies the definition of the shortest
biconnection. Suppose that there is a unique shortest path $P=P_1$.
Then $P_2$ is a unique shortest path between $u$ and $v$ in 
$G\setminus(V(P_1)\setminus\{u,v\})$ for else $C$ would not be a unique
shortest cycle via $u$ and $v$. Therefore $C$ is the shortest biconnection
also in this case.

Assume now that $C$ is boundary-like and that $u$ and $v$ are as above.
We now have to show that $C$ is the shortest cycle via $u$ and $v$.
If $|P_1|=|P_2|$, this easily follows from the fact that $C$ is the
shortest biconnection of $u$ and $v$. 

Suppose $|P_1|<|P_2|$.
Assume that $C'$ is another cycle via $u$ and $v$ with $|C'|\le|C|$.
Denote the paths into which $u$ and $v$ split
$C'$ by $P'_1$ and $P'_2$. As easily seen, neither of them  coincides with $P_1$.
Hence both $|P'_1|$ and $|P'_2|$ exceed $|P_1|$ and both do not exceed $|P_2|$.
By the definition of a biconnection, both $P'_1$ and $P'_2$ have, except $u$ and $v$,
other common vertices with $P_1$. Let $u',v'\notin\{u,v\}$ be two vertices
in $P_1$ such that $u'\in P'_1$, $v'\in P'_2$, and no vertex between $u'$ and
$v'$ belongs to $P'_1$ or $P'_2$.

Let $u'$ and $v'$ split $C$ into paths $Q_1$ and $Q_2$ where 
$Q_1$ is a part of $P_1$ and hence $|Q_1|=d(u',v')$.
Let $u'$ and $v'$ split $C'$ into paths $Q'_1$ and $Q'_2$.
Since $|Q'_i|\ge d(u',v')$ for both $i$, for both $i$ we have
$|Q'_i|=|C'|-|Q'_{3-i}|\le|C|-|Q_1|=|Q_2|$. 
Assume that $u'$ and $v'$ are non-adjacent.
Coupled with $Q_1$,
either $Q'_i$ shows that $C$ cannot be the shortest biconnection of
$u'$ and $v'$, a contradiction.

Assume that $u'$ and $v'$ are adjacent. Without loss of generality suppose that
$u'$ is between $u$ and $v'$ in $P_1$. Consider the part of $P'_1$ from
$u'$ to $v$. Let $v''$ be the nearest to $u'$ vertex on this arc which
belongs to $P_1$. It is possible that $v''=v$ but not $v''=v'$ because $v'\in P'_2$.
Thus $d(u',v'')\ge2$. Denote the arc of $P'_1$ between $u'$ and $v''$ by $R'$.
We have $|R'|\le|P'_1|\le|P_2|$ which is strictly less than the length of the 
longer arc of $C$ between $u'$ and $v''$. Thus, coupled with the arc of $P_1$
between $u'$ and $v''$ the $R'$ shows that $C$ cannot be the shortest biconnection
of between $u'$ and $v''$, a contradiction which completes the proof.
\end{proof}

\begin{definition}[a pseudo-BOP graph]\label{def:psBOP}\rm
Call $G$ a \emph{pseudo-BOP graph} if the following conditions are met.
\begin{enumerate}
\item
$G$ has no isolated vertex.
\item
Every edge $e\in E(G)$ belongs to at least one but at most two pseudo-facial cycles.
If $e$ belongs to exactly one, we will call it \emph{outer}.
\item
If edges $\{u,v\}$ and $\{v,w\}$ lie on pseudo-facial cycles and both are
outer, then the vertex $v$ belongs to no other pseudo-facial cycle.
\end{enumerate}

Let $G$ be a pseudo-BOP graph. In a way coherent with Definition \ref{def:bop}, 
let $\oo G$ denote the set of outer edges of~$G$.
\end{definition}

Note that every BOP graph is pseudo-BOP. Other examples of pseudo-BOP graphs
are all plane graphs with every face, maybe except the unbounded one,
being a triangle, and the modification of $P_n\times P_2$ to a Moebius strip.

\section{Distinguishing a BOP graph from a non-pseudo-BOP graph}\label{s:bopvsnon}

Our goal in this section is to prove the following technical result.

\begin{lemma}[Main Lemma 2A]\label{lem:ml1a}
If $G$ is a BOP graph and $G'$ is not a pseudo-BOP graph, then
Spoiler wins $\ehr G{G'}$ in less than $2\,\log f(G)+2\,\log r(\dual G)+9$ moves.
\end{lemma}

We first extend our collection of Spoiler's tricks in $\ehr G{G'}$ 
from Section \ref{s:ehr}. In a series of lemmas below we assume that
$u,v\in V(G)$, $u',v'\in V(G')$ and each pair
$u,u'$ and $v,v'$ is selected in a round.

\begin{lemma}[Metric Threat 2]\label{lem:mt2}
Suppose that $d(u,v)=d(u',v')$.
Then, whenever Spoiler selects a vertex $w$ on a shortest path between
$u$ and $v$, Duplicator in response selects $w'$ on a shortest path between
$u'$ and $v'$ with $d(w',u')=d(w,u)$ and $d(w',v')=d(w,v)$ 
or loses in less than $\lceil\log d(u,v)\rceil$ moves.
\end{lemma}

Lemma \ref{lem:mt2} is a direct consequence of Lemma \ref{lem:mt1} (Metric Threat~1).

\begin{lemma}[Metric Threat 3]\label{lem:mt3}
Suppose that $u$ and $v$ have a unique shortest biconnection $C$.
Then $u'$ and $v'$ have a unique shortest biconnection $C'$
with $|P'_1|=|P_1|$ and $|P'_2|=|P_2|$ or otherwise
Spoiler wins in less than $2\log|C|+1$ moves.
\end{lemma}

\begin{proof}
We suppose that $d(u,v)=d(u',v')$; Otherwise we have Metric Threat~1.

\case 1{$|P_1|=|P_2|$.} 
We split it in three subcases.

\subcase{1.1}{There are two shortest paths between $u'$ and $v'$
that have a common vertex $w'\notin\{u',v'\}$.}
Assume that the paths diverge at $w'$ and go further through
distinct vertices $z'$ and $z''$. Spoiler first selects $w'$ and
Duplicator should respond with a vertex $w$ at the same distances
from $u$ and $v$ on $P_1$ or on $P_2$. Then Spoiler selects $z'$ and $z''$.
Irrespectively of which neighbors of $w$ are selected of Duplicator,
she faces Metric Threat 2 and loses in less than $\log|P_1|+1$ moves.

\subcase{1.2}{There are three non-crossing shortest paths between $u'$ and $v'$.}
Spoiler selects three vertices which are the neighbors of $u'$ in these paths.
This poses Metric Threat~2.

\subcase{1.3}{There is a unique shortest path between $u'$ and $v'$.}
Spoiler selects the neighbors of $u$ in $P_1$ and $P_2$ thereby posing
Metric Threat~2.

\case 2{$|P_1|<|P_2|$.}

\subcase{2.1}{There are at least two shortest paths $P'$ and $P''$ 
between $u'$ and $v'$.}
There are two vertices $w'\in P'$ and $w''\in P''$ at the same distances
between $u'$ and $v'$. Spoiler selects $w'$ and $w''$ thereby posing
Metric Threat~2.

\subcase{2.2}{There is a unique shortest path $P'_1$ between $u'$ and $v'$.}
Spoiler forces play on $H=G\setminus(V(P_1)\setminus\{u,v\})$ and
$H'=G'\setminus(V(P'_1)\setminus\{u',v'\})$, that is, from now on
he never moves inside $P_1$ or $P'_1$. Once Duplicator moves inside
$P_1$ or $P'_1$, she faces Metric Threat 2. We hence suppose
that she never does so. If $d_{H'}(u',v')\ne d_{H}(u,v)=|P_2|$ or if
in $H'$ there are two shortest paths between $u'$ and $v'$,
Spoiler wins exactly as above. The game lasts less than
$\log|P_2|+\log|P_1|+3\le\log|C|+1$ rounds (the term $\log|P_1|$
appears because Duplicator, before her loss in $\ehr H{H'}$,
can start moving somewhere in $P_1$ or $P'_1$).
\end{proof}

\begin{lemma}[Metric Threat 4]
Given a vertex $w$ in $C$ which is the shortest biconnection of $u$ and $v$,
we define the \emph{coordinates of $w$} to be the pair $(i,d(w,u))$
where $i$ is as follows: If $|P_1|=|P_2|$, then $i$ is the empty word;
Otherwise $i$ is determined by the condition $w\in P_i$.

Suppose that $C$ is the shortest biconnection of $u$ and $v$ in $G$,
$C'$ is the shortest biconnection of $u'$ and $v'$ in $G'$,
and $|P_i|=|P'_i|$, $i=1,2$. Furthermore, suppose that $w\in V(G)$
and $w'\in V(G')$ are selected in a round and that $w\in C$.
Then $w'\in C'$ and has the same coordinates as $w$
or otherwise Spoiler wins in less than $2\log|C|+1$ moves.
\end{lemma}

\begin{proof}
Suppose that the coordinates of $w'$ in $C'$ differ from the
the coordinates of $w$ in $C$. If $|P_1|=|P_2|$, this is Metric Threat 2.
Let $|P_1|<|P_2|$. If $w\in P_1$, this is Metric Threat 1.
If $w\in P_2$, this is Metric Threat 1 in the graphs $H$ and $H'$
as in the proof of Lemma \ref{lem:mt3} and Spoiler wins as explained there.
\end{proof}

\begin{lemma}[Metric Threat 5]
Suppose that $u$ and $v$ are non-adjacent and lie on a pseudo-facial cycle $C$.
Then $u'$ and $v'$ lie on a pseudo-facial cycle $C'$
with $|P'_1|=|P_1|$ and $|P'_2|=|P_2|$ or otherwise
Spoiler wins in less than $2\log|C|+3$ moves.
\end{lemma}

\begin{proof}
Recall that, by Lemma \ref{lem:pfeqbl}, the notions of a pseudo-facial
and a boundary-like cycle coincide.
Unless we have Metric Threat 3, $u'$ and $v'$ lie on the shortest biconnection 
$C'$ with $|P'_1|=|P_1|$ and $|P'_2|=|P_2|$. Suppose that the cycle $C'$
is not pseudo-facial. This means that there are two non-adjacent 
vertices $s',t'\in C'$ for which $C'$ is not the shortest biconnection.
Let $s$ and $t$ be the vertices with the same coordinates on $C$.
Spoiler selects $s$ and $t$ and Duplicator should respond with $s'$ and $t'$
for else she faces Metric Threat 4. Once Duplicator does so, she
faces Metric Threat~3.
\end{proof}

\begin{lemma}\label{lem:ffggempt}
If $\ff G\ne\emptyset$ and $\ff{G'}=\emptyset$, then Spoiler wins $\ehr{G}{G'}$
in less than $2\log f(G)+5$ moves. 
\end{lemma}

\begin{proof}
If $G$ contains a triangle Spoiler selects its vertices and wins.
Otherwise Spoiler selects two non-adjacent vertices in a pseudo-facial cycle
and wins making Metric Threat~5.
\end{proof}

\begin{lemma}\label{lem:maxfgfg}
If $G$ is pseudo-BOP, $G'$ is not, and $\ff{G'}\ne\emptyset$, then Spoiler 
wins $\ehr{G}{G'}$ in $2\log\max\{f(G),f(G')\}+8$ moves.
\end{lemma}

\begin{proof}
We trace through the items of Definition \ref{def:psBOP}.

\case 1{$G'$ has an isolated vertex.}
Spoiler wins in 2 moves.

\case 2{$G'$ has an edge $e'$ that belongs to no pseudo-facial cycle.}
Spoiler selects $e'=\{t',v'\}$. Let Duplicator respond with
the $e=\{t,v\}$, an edge of $G$. The $e$ belongs to a pseudo-facial
cycle $C$. Spoiler selects $u$, the other neighbor of $t$ in $C$.
If $u$ and $v$ are adjacent, this is Duplicator's loss. Otherwise
Duplicator faces Metric Threat 5 and loses in less than extra 
$2\log|C|+3$ moves.

\case 3{$G'$ has an edge $e'$ that belongs to three pseudo-facial cycles.}
The first two rounds are as in the preceding case. In the next three
rounds Spoiler selects $u'_1$, $u'_2$, and $u'_3$, the neighbors of $t'$
different from $v'$ in the three pseudo-facial cycles. By Proposition
\ref{prop:atmosttwo} these vertices are pairwise distinct. Whatever
Duplicator responds, one of the pairs $u'_i,v$ poses Metric Threat~5
and Spoiler wins in less than $2\log f(G')+3$ extra moves.

\case 4{$G'$ has outer edges $\{u',v'\}$ and $\{v',w'\}$ lying on 
pseudo-facial cycles such that the vertex $v'$ belongs to some other 
pseudo-facial cycle $C'$.}
Spoiler selects $u'$, $v'$, and $w'$. Let Duplicator respond with
$u$, $v$, and $w$. If at least one of the edges $\{u,v\}$ and $\{v,w\}$
is not outer, Spoiler wins in less than $2\log f(G)+5$ extra moves similarly 
to Case 3.
Otherwise $v$ belongs to no other pseudo-facial cycle.
Spoiler selects a vertex $z'\in C'$ non-adjacent to $v'$
(if $C'$ is a triangle, the life is much easier).
Denote Duplicator's response by $z$. If $z$ belongs to the pseudo-facial
cycle containing $\{u,v\}$ or to the pseudo-facial
cycle containing $\{v,w\}$, in the next move Spoiler is able to create
Metric Threat 4 or 5 (because $u',w'\notin C'$) and win in
less than $2\log f(G)+3$ next moves. Otherwise $z$ belongs to no 
pseudo-facial cycle which poses Metric Threat 5 and hence Spoiler wins
in less than $2\log |C'|+3$ next moves.
\end{proof}

\begin{lemma}\label{lem:fgnefg}
Let $G$ be a BOP graph and $G'$ be an arbitrary graph with $\ff{G'}\ne\emptyset$.
Suppose that $f(G)\ne f(G')$. Then Spoiler wins $\ehr{G,G'}$ in less than
$2\log f(G)+2\log r(\dual G)+9$ moves.
\end{lemma}

\begin{proof}
Set $f=f(G)$ and $r=r(\dual G)$. Suppose that $G'$ has no isolated vertex and no leaf.
Given a vertex $v$ of $G$ or $G'$, we define an \emph{adjoining relation}
on $\Gamma(v)$: We say that $a,b\in\Gamma(v)$ \emph{adjoin} if the edges
$\{v,a\}$ and $\{v,b\}$ lie on the same pseudo-facial cycle of length
no more than $f$. In the BOP graph $G$ the adjoining relation induces 
a path through all $\Gamma(v)$ whose endpoint are the neighbors of $v$
in the outer facial cycle of~$G$.

\begin{claim}\label{cl:adjoin}
For every $v'\in V(G')$, the adjoining relation induces
a path through all $\Gamma(v')$ or otherwise Spoiler wins
in less than $2\log f+9$ moves.
\end{claim}

\begin{subproof}
Suppose that for some $v'\in V(G')$ the first part of the claim is not true.
In the first round Spoiler selects $v'$. Let $v\in V(G)$ denote Duplicator's
response.

\case 1{$\Gamma(v')$ contains four vertices $a',b'_1,b'_2,b'_3$
such that $a'$ adjoins every $b'_i$.}
Then $G'$ is not a pseudo-BOP graph which is witnessed by the edge
$\{v',a'\}$ belonging to three pseudo-facial cycles, each of length at most
$f$. As in the proof of Lemma \ref{lem:maxfgfg}, Case 3,
Spoiler wins in less than $2\log f+8$ moves at total.

\case 2{$\Gamma(v')$ contains a vertex $a'$ adjoining no other vertex 
in $\Gamma(v')$.}
Spoiler selects $a'$. Let Duplicator respond with 
$a\in\Gamma(v)$. In the next move Spoiler selects $b\in\Gamma(v)$
adjoining $a$, which poses Metric Threat 5 irrespective of
Duplicator's response in $\Gamma(v')$. Implementing it, Spoiler
wins in less than $2\log f+3$ extra moves.

\case 3{Every vertex in $\Gamma(v')$ has two adjoining vertices in $\Gamma(v')$.}
Spoiler selects a vertex $a\in\Gamma(v)$ with exactly one adjoining vertex.
Denote Duplicator's response in $\Gamma(v')$ by $a'$. Then Spoiler selects
two vertices $b'_1,b'_2\in\Gamma(v')$ both adjoining $a'$. Let Duplicator
respond with $b_1,b_2\in\Gamma(v)$. For $i=1$ or $i=2$, $a'$ and $b'_i$ lie
on a pseudo-facial cycle of length at most $f$ while this is not so for $a$ and $b_i$.
This is Metric Threat 5 and Spoiler
wins in less than $2\log f+3$ extra moves.

\case 4{$\Gamma(v')$ contains three vertices $a'_1,a'_2,a'_3$
with exactly one adjoining vertex each.}
Spoiler selects these three. One of Duplicator's responses, $a_i$,
has two adjoining vertices $b_1,b_2\in\Gamma(v)$ while $a'_i$ has only one.
Spoiler selects $b_1$ and $b_2$. For one of Duplicator's responses 
$b'_j\in\Gamma(v')$, $a'_i$ and $b'_j$ are not adjoining while $a_i$ and $b_j$ are.
This is again Metric Threat 5 and Spoiler
wins in less than $2\log f+3$ extra moves.
\end{subproof}

In the sequel we suppose that the first part of Claim \ref{cl:adjoin}
is true.

The case of $f(G') < f(G)$ is simple. Spoiler selects two non-adjacent
vertices on the longest pseudo-facial cycle in $G$ and wins in less than
$2\log f+3$ next moves implementing Metric Threat 5. We hence suppose
that $f(G') > f(G)$. The same trick does not apply as we wish that our bound
does not depend on~$G'$.

Fix $C'\in\ff{C'}$ with $|C'|>f$. Let $u'$, $v'$, and $w'$ be three successive
vertices on $C'$ and $d=\deg v'$. Going through pseudo-facial cycles
in the order given by the adjoining relation on $\Gamma(v')$, one can
come from $u'$ to $w'$ not visiting $v'$ via a route of length at most
$(d-1)(f-2)$. It follows that $|C'|<df$, which allows Spoiler to win
as above in less than $2\log|C'|+5<2\log d+2\log f+5$ moves.
If $C'$ contains a vertex $v'$ of degree at most $r+1$ in $G'$, we are done.

In what follows we suppose that $C'$ consists of vertices of degree at least $r+2$.

With each $v$ in $G$ or $G'$ we associate two number sequences, $s_1(v)$
and $s_2(v)$. Order the bunch of pseudo-facial cycles of length at most $f$
going through the vertex $v$ in accordance with the adjoining relation.
Then $s_1(v)$ and $s_2(v)$ are sequences of the sizes of cycles in the bunch
in the direct and the reverse order. Define $S(v)$ to be the set of all
subwords of $s_1(v)$ and $s_2(v)$ having length $r$ (a subword is a subsequence
of consecutive elements).

\begin{claim}\label{cl:svsv}
Suppose that $v$ and $v'$ are selected in a round and $S(v)\ne S(v')$.
Then Spoiler wins in at most $2\log f+\log r+6$ moves.
\end{claim}

\begin{subproof}
Without loss of generality suppose that there is an 
$\alpha\in S(v)\setminus S(v')$. Consider the path of the adjoining relation
corresponding to $\alpha$ and denote its endvertices by $u$ and $w$.
Spoiler selects these two vertices. Let Duplicator respond with
$u',w'\in\Gamma(v')$. Consider the path of the adjoining relation
between $u'$ and $w'$ and denote the corresponding sequence of the sizes
of pseudo-facial cycles of length at most $f$ by $\alpha'$.
If the length of $\alpha'$ is not equal to $r$, the length of $\alpha$,
then Spoiler plays in $\Gamma(v)\cup\Gamma(v')$ applying the halving
strategy distinguishing between two paths of different length, namely,
the adjoining-relation paths between $u$ and $w$ and between $u'$ and $w'$.
In less than $\log(r+1)+1$ moves he forces appearance of $a,b\in\Gamma(v)$
and $a',b'\in\Gamma(v')$ such that $a$ and $b$ are adjoining iff $a'$ and $b'$
are not. If one of the pairs is adjacent, this is an immediate win
of Spoiler. Otherwise, this is Metric Threat 5 and Spoiler needs
less than $2\log f+3$ extra moves to win. 

Suppose therefore that $\alpha$ and $\alpha'$ have the same length.
Let them differ at the $l$-th position. Denote the corresponding entries
by $\alpha_l$ and $\alpha'_l$. Let $F$ be the corresponding pseudo-facial
cycle of length $\alpha_l$ and $F'$ the corresponding pseudo-facial
cycle of length $\alpha'_l$. Denote the neighbors of $v$ in $F$ by $s$ and $t$
and the neighbors of $v'$ in $F'$ by $s'$ and $t'$. Spoiler selects $s$ and $t$.
If Duplicator does not respond with $s'$ and $t'$, she loses in less than
$(\log l+1)+(2\log f+3)$ moves as above. If Duplicator responds with $s'$ and $t'$,
she loses in less than $2\log f+3$ moves needed for Spoiler to implement
Metric Threat~5.
\end{subproof}

In view of Claim \ref{cl:svsv}, we will assume that for every $v'\in C'$
there is a $v\in V(G)$ with $S(v)=S(v')$. We will call such a $v$
a \emph{friend of the $v'$}. By the definition of $r$-fineness, for every
two vertices $v_1$ and $v_2$ in $G$ each of degree at least $r+2$, we have
$S(v_1)\ne S(v_2)$. This implies that a vertex $v'\in C'$ has at most one
friend of degree at least $r+2$. Such a friend will be referred to as
the \emph{big friend of~$v'$}.

\begin{claim}\label{cl:notbig} 
Suppose that $v\in G$ and $v'\in C'$ are selected in a round and that
$v$ is not a big friend of $v'$. Then Spoiler wins in 
at most $2\log f+\log r+6$ moves.
\end{claim}

\begin{subproof}
If $S(v)\ne S(v')$, Spoiler wins by Claim \ref{cl:svsv}. Otherwise $\deg v\le r+1$,
while $\deg v'\ge r+2$. In this case Spoiler wins in at most $2\log f+\log r+6$ moves
using the adjoining-relation paths on $\Gamma(v)$ and $\Gamma(v)$
that have different length.
\end{subproof}

On the account of Claim \ref{cl:notbig}, we will suppose that
every vertex $v'$ of $C'$ has a big friend $v$ in $G$ (which is unique) and that, 
whenever Spoiler selects $v'$, Duplicator responds with the $v$.
Under this assumption, if two different vertices in $C'$ have the same
big friend, Spoiler wins in 2 moves. We therefore will also suppose that
the big friendship is a one-to-one correspondence between $V(C')$
and a set $B\subset V(G)$.

Assume for a while that the big friendship
is a partial isomorphism between $G$ and $G'$.
Denote the cycle that $B$ spans in $G$ by $C$. 
As $C$ is chordless, it is a facial cycle in $G$ contradicting the inequality
$|C|=|C'|>f(G)$.

Thus, $C'$ must contain two vertices $u'$ and $v'$ 
such that their big friends $u$ and $v$ have different adjacency.
Spoiler wins selecting $u'$ and $v'$ unless Duplicator breaks the big friendship.
\end{proof}

\begin{proofof}{Main Lemma 2A (Lemma \ref{lem:ml1a})}
If $\ff{G'}=\emptyset$, then Spoiler has a fast win by Lemma \ref{lem:ffggempt}.
Let $\ff{G'}\ne\emptyset$. Then Spoiler has a fast win by Lemma \ref{lem:fgnefg} 
if $f(G')\ne f(G)$ and by Lemma \ref{lem:maxfgfg} otherwise. 
\end{proofof}

\section{Distinguishing a BOP graph from a non-isomorphic pseudo-BOP graph}%
\label{s:bopvspseudo}

We here proceed with the second ingredient of the proof of 
Main Lemma 2 (Lemma \ref{lem:ml1}),
Lemma \ref{lem:ml1b} below. In the preceding section we were able to distinguish
a BOP graph $G$ from a non-pseudo-BOP graph $G'$. Now we need to distinguish $G$
from a non-isomorphic pseudo-BOP graph $G'$ by proving Inequality \ref{eq:ggfacing}.
The first what we have to do is to extend the definition of $\facing$ to pseudo-BOP
graphs.

\begin{definition}[the facing structure for pseudo-BOP graphs]\rm
Let $G$ be a pseudo-BOP graph. We define a structure $\facing G=\langle H,L\rangle$,
which is a graph $H$ with layout $L$, literally as in Definition \ref{def:facing},
where the notation $\ff G$ and $\oo G$ is extended to pseudo-BOP graphs.
As in Definition \ref{def:facing}, we simultaneously define the \emph{crossing relation}.
We also keep the notation $\dual G$ for~$H$.
\end{definition}

\begin{proposition}\label{prop:iftree}
Let $G$ be a pseudo-BOP graph whose $\dual G$ is a tree. Then $G$ is BOP.
\end{proposition}

\begin{proof}
We proceed by induction on $k=|\ff G|$. In the base case of $k=1$
we easily see that $G$ is a cycle. Let $k\ge2$. Fix a $C\in\ff G$
that is a leaf in $\dual G$. The latter condition means that all
edges of $C$ except one, say $\{u,v\}$, are outer. By the definition
of a pseudo-BOP graph, every vertex in $V(C)\setminus\{u,v\}$ has
degree 2 in $G$ while $\deg u,\deg v\ge4$. Let $H=G\setminus(V(C)\setminus\{u,v\})$.
It is not hard to see that $\ff H=\ff G\setminus\{C\}$. On the account
of this fact, we easily see that $H$ is a pseudo-BOP graph in which
$\{u,v\}$ is outer. By the induction assumption $H$ is BOP.
Since glueing a new cycle to a BOP graph at an outer edge gives
a BOP graph, $G$ is BOP.
\end{proof}

\begin{proposition}\label{prop:isoiff}
Let $G$ be a BOP graph and $G'$ a pseudo-BOP graph. Then $G\cong G'$
if and only if $\facing G\cong\facing{G'}$, where the latter isomorphism is between
2-relation structures.
\end{proposition}

\begin{proof}
By Propositions \ref{prop:facingBOP} and~\ref{prop:iftree}.
\end{proof}

Note that the statement below needs Proposition \ref{prop:isoiff} to be sound.

\begin{lemma}[Main Lemma 2B]\label{lem:ml1b}
Let $G$ be a BOP graph and $G'$ be a non-isomorphic pseudo-BOP graph. Then
$$
D(G,G') < 3\,D(\facing G,\facing{G'})+2\,\log f(G)+2\,\log r(\dual G)+5.
$$
\end{lemma}

The proof is very short modulo two lemmas.
If $f(G')>f(G)$, then we are done by Lemma \ref{lem:fgnefg}.
If $f(G')\le f(G)$, then Main Lemma 2B immediately follows from
another lemma which we state and prove below. 

\begin{lemma}
Let $G$ be a BOP graph and $G'$ be a non-isomorphic pseudo-BOP graph. Then
$$
D(G,G') < 3\,D(\facing G,\facing{G'})+2\,\log\max\{f(G),f(G')\}+5.
$$
\end{lemma}

\begin{proof}
We have to design a winning strategy for Spoiler in $\ehr G{G'}$.
Spoiler will simulate $\ehr{\facing G}{\facing G'}$.
We will refer to $\ehr G{G'}$ as \emph{the game} and to
$\ehr{\facing G}{\facing G'}$ as \emph{the protogame}.
Let $\facing G=\langle H,L\rangle$ and $\facing G'=\langle H',L'\rangle$.

The simulation goes as follows.
\begin{itemize}
\item
Whenever Spoiler selects a vertex $C\in\ff G\cup\ff{G'}$
in the protogame, in the game he selects all the vertices of $C$
if it is a triangle and only 2 non-adjacent vertices of $C$ otherwise.
Note that this unambiguously specifies the cycle $C$ by Proposition 
\ref{prop:atmosttwo}.
\item
Whenever Spoiler selects a vertex $e\in\oo G\cup\oo{G'}$
in the protogame, in the game he selects the endvertices of~$e$.
\end{itemize}
Note now that, if Duplicator's play does not correspond to the protogame
in this manner, Spoiler has a fast win.
\begin{itemize}
\item
If Duplicator responds not with a triangle or not with two
non-adjacent vertices lying on a pseudo-facial cycle, she
loses immediately or faces
Metric Threat 5 and loses in less than $2\log\max\{f(G),f(G')\}+3$ moves.
\item
If Duplicator responds with a non-edge or with a non-outer edge $e'$,
then she loses immediately or Spoiler wins in less than 
$2\log\max\{f(G),f(G')\}+5$ moves exhibiting two pseudo-facial cycles
containing $e'$ (like to Case 2 or 3 in the proof of Lemma \ref{lem:maxfgfg}).
\end{itemize}

Suppose that Duplicator always bewares of these threats. 
Then a round of the protogame corresponds to two or three rounds of
the game. It will be convenient to have this correspondence more uniform
and we assume that, whenever Spoiler has to make two moves, he makes three
just by repeating the second move once again. Thus, $3k$ rounds of the
game correspond to $k$ rounds of the protogame.
After completing $3D(\facing G,\facing G')$ rounds (when the protogame
is won by Spoiler) the game comes in the phase of endgame.
Several endgames are possible depending on the terminal position of
the protogame.

\fingame 1{A partial isomorphism between $H$ and $H'$ is broken.}
We split our analysis in two cases.

\fingame{1.1}{$C,D\in\ff G$ and $C',D'\in\ff{G'}$ are selected in the protogame
and have different adjacency in $H$ and $H'$ respectively.}
Without loss of generality assume that $C$ and $D$ share an edge $e$
while $C'$ and $D'$ have no edge in common. Spoiler selects the endvertices
of $e$. At least one of the vertices selected in response by Duplicator
will not belong to $C'$ or $D'$. This is either immediate Duplicator's loss
or at least Metric Threat 4 and loss in less than $2\log\max\{f(G),f(G')\}+1$
moves.

\fingame{1.2}{$C\in\ff G$, $e\in\oo G$ and $C'\in\ff{G'}$, $e'\in\oo{G'}$ 
are selected in the protogame
and have different adjacency in $H$ and $H'$ respectively.}
Without loss of generality assume that $C$ contains $e$ while $C'$
does not contain $e'$.
As at least one of the endvertices of $e'$ does not lie on $C'$,
Duplicator again loses in less than $2\log\max\{f(G),f(G')\}+1$ moves.

To analyze the endgame of a different type, we use the following
observation.

\begin{claim}\label{cl:cross}
Suppose that in the protogame adjacent vertices $A,B\in\ff G\cup\oo G$
and $A',B'\in\ff{G'}\cup\oo{G'}$ have been selected.
Let Spoiler select the endvertices of the crossing edge $\{A,B\}^*$ in the game.
Then Duplicator should respond with the endvertices of the crossing edge 
$\{A',B'\}^*$ for else she loses in less than $2\log\max\{f(G),f(G')\}+1$ moves.
A symmetric claim with $G$ and $G'$ interchanged is also true.
\end{claim}

\begin{subproof}
\case 1{$A,B\in\ff GG$.}
Note that $\{A,B\}^*$ (resp.\ $\{A',B'\}^*$) is a unique common edge of the cycles 
$A$ and $B$ in $G$ (resp.\ $A'$ and $B'$ in $G'$). Assume that,
in response to $\{A,B\}^*$, Duplicator selects the endvertices of an edge
different from $\{A',B'\}^*$. At least one of the two vertices selected by her
does not belong to $A'$ or $B'$ while the corresponding vertex in $G$
belongs to both $A$ or $B$. Duplicator loses either immediately or faces 
Metric Threat 4 and loses in less than $2\log\max\{f(G),f(G')\}+1$ moves.

\case 2{$A\in\ff G$, $B=e\in\oo G$, $A'\in\ff{G'}$, $B'=e'\in\oo{G'}$.} 
Note that the edge $\{A,e\}^*=e$ belongs to the cycle $A$ and the edge $\{A',e'\}^*=e'$
belongs to the cycle $A'$.
Selection of $e$ and $e'$ in the protogame means that the endvertices
of $e$ and $e'$ are selected in the game. Now we suppose that
Spoiler selects the endvertices of $\{A,e\}^*=e$ (once again)
but Duplicator does not respond with the endvertices of $\{A',e'\}^*=e'$.
Therefore the equality relation is violated and Duplicator has lost.
\end{subproof}

\fingame 2{A partial isomorphism between $L$ and $L'$ is broken.}
Suppose this is witnessed by selection of vertices $U,V\in\ff G\cup\oo G$
and $U',V'\in\ff{G'}\cup\oo{G'}$ in the protogame.
Without loss of generality we assume that $U$ and $V$ are in relation $L$
while $U'$ and $V'$ are not in relation $L'$. The opposite case is treated
by the symmetric argument: We need the fact that $G$ is BOP (not just pseudo-BOP)
only to have $\facing G\not\cong\facing G'$ and do not use this fact any more.

\fingame{2.1}{$U,V\in\Gamma(W)$ for some $W$ in $H$.}
Spoiler selects $W$ in the protogame. 
Denote Duplicator's response by $W'$.
Suppose that $U',V'\in\Gamma_{H'}(W')$ because otherwise we arrive at Endgame 1.
By assumption and the definition of a layout, 
the edges $\{W,U\}^*$ and $\{W,V\}^*$ are adjacent in $H$
while $\{W',U'\}^*$ and $\{W',V'\}^*$ are non-adjacent in $H'$.
Spoiler selects the three endvertices of $\{W,U\}^*$ and $\{W,V\}^*$.
Duplicator can respond with at most three of the four endvertices of
$\{W',U'\}^*$ and $\{W',V'\}^*$. By Claim \ref{cl:cross} she loses
in less than $2\log\max\{f(G),f(G')\}+1$ next moves.

\fingame{2.2}{$U,V\in\Gamma(W)$ for no $W$ in $H$.}
By the definition of a layout, between $U$ and $V$ in $H$ there is a path
of length 3. We suppose that such a path exists in $H'$ between $U'$ and $V'$
because otherwise Spoiler is able to force Endgame 1 in two moves.
Let $A,B\in\ff G$ and $A',B'\in\ff{G'}$ be
the intermediate vertices in these paths.
By assumption and the definition of a layout, 
the edges $\{U,A\}^*$ and $\{B,V\}^*$ are adjacent in $H$
while $\{U',A'\}^*$ and $\{B',V'\}^*$ are non-adjacent in $H'$.
Spoiler selects the three endvertices of $\{U,A\}^*$ and $\{B,V\}^*$.
Duplicator can respond with at most three of the four endvertices of
$\{U',A'\}^*$ and $\{B',V'\}^*$. By Claim \ref{cl:cross} she loses
in less than $2\log\max\{f(G),f(G')\}+1$ next moves.
\end{proof}

We are now able to complete the task of this section.

\begin{proofof}{Main Lemma 2 (Lemma~\ref{lem:ml1})}
Let $G$ be a BOP graph.
Given a graph $G'$ non-isomorphic to $G$, we have to estimate $D(G,G')$
from above.
If $G'$ is not pseudo-BOP, then 
$$
D(G,G')<2\,\log f(G)+2\,\log r(\dual G)+9
$$
by Main Lemma 2A (Lemma \ref{lem:ml1a}).
If $G'$ is pseudo-BOP, then Main Lemma 2B (Lemma \ref{lem:ml1b}) gives
\begin{eqnarray*}
D(G,G')&<&3\,D(\facing G,\facing{G'})+2\,\log f(G)+2\,\log r(\dual G)+5\\
&\le&3\,D(\facing G)+2\,\log f(G)+2\,\log r(\dual G)+5.
\end{eqnarray*}
Since $D(G)=\max_{G'}D(G,G')$, the claimed bound follows.
\end{proofof}

\section{Proof of Theorem \protect\ref{thm:maiin}}\label{s:uppeer}

Main Lemmas 1 and 2 (Lemmas \ref{lem:ml2} and \ref{lem:ml1} respectively) 
imply that for every BOP graph $G$ we have
\begin{eqnarray*}
D(G)&<&3\,D(\facing G)+2\,\log f(G)+2\,\log r(\dual G)+5\\
&\le&11\,\log r(\dual G)+5\,\log\Delta(\dual)+59,
\end{eqnarray*}
as required
(note that $\Delta(\dual G)=f(G)$).

\section{Future work}\label{s:futur}

Let $G_n$ be a random biconnected outerplanar graph on $n$ labeled vertices.
In other words, we fix a convex cycle of length $n$ in a plane, label its
vertices at random, and add some non-crossing chords so that
all sets of non-crossing chords are equiprobable.
We conjecture that $D(G_n)=\Theta(\log\log n)$ with probability $1-o(1)$
and believe that the upper bound here can be derived from Theorem \ref{thm:maiin}.
This is so if the following two combinatorial hypotheses are true. First,
$\Delta(\dual G_n)=O(\log^{O(1)}n)$ with probability $1-o(1)$.
Second, $r(\dual G_n)=O(\log^{O(1)}n)$ with probability $1-o(1)$.
The double-logarithmic lower bound would follow from the hypothesis that, with probability
$1-o(1)$, $G_n$ has an induced path of length $\log^{\Omega(1)}n$ whose vertices
have degree 2 in the graph.

\end{document}